\documentstyle[12pt]{article}
\titlepage

\title{The Geometry of a pair of second-order ODEs and Euclidean spaces}
\author{Richard Atkins} 
\date{}
\newtheorem{fact}{Fact}

\newtheorem{theorem}[fact]{Theorem}

\begin{document}
\maketitle 
Department of Mathematics, Trinity Western University, Glover Road, Langley, 
British Columbia V2Y 1Y1, Canada \\ 
\maketitle
\newpage
\section{Introduction}
The geometry of a manifold is given by a metric, which defines a notion of 
distance between two points. Paths of shortest length connecting points are 
obtained as the critical curves of the functional variation of the integral
defining arclength. Functional variation of this integral yields Euler-Lagrange
equations which are a system of ordinary differential equations of second order, 
whose solutions are the geodesics. Thus associated with geometry is a system of 
ODEs. This paper seeks to answer the inverse problem: when does a system of ODEs 
represent the paths of shortest length of a metric? That is, we wish to establish
when ordinary differential equations exhibit an underlying geometry. We shall
not be so ambitious as to attempt a solution on manifolds of arbitrary dimension
and endowed with a general metric but shall restrict ourselves to the case of a pair
of second-order ODEs on a (two-dimensional) surface and ask when the underlying geometry
is flat, that is, a Euclidean space. We are concerned only with the solutions of the 
ODEs up to reparameterization since they serve merely to describe the paths of shortest
length on the surface. Geodesics however, are not invariant with respect to general
changes of parameter, so it shall be necessary to incorporate reparameterization
in the precise definition of the problem. The formulation of the equivalence problem
is the content of sections 1 and 2. In section 3, Cartan's method of equivalence is
employed up to the level of the first normalization for generic ODEs. In section 4,
we obtain the solution in the Euclidean case. The equivalence method is carried through
in section 5, for generic pairs of second-order ODEs with the result that the symmetry 
of the system produces 24 invariant functions. 

Before proceeding, it is intructive to consider the simpler  problem without considerations
of parameterization: when are the solutions
to a pair of second-order ordinary differential equations
\[ d^{2}y^{1}/d^{2}t = f(y,\dot{y},t)  \hspace*{.3in} and \hspace*{.3in}
d^{2}y^{2}/d^{2}t = g(y,\dot{y},t),  \hspace{1in} (1) \]
locally, the geodesics of some Euclidean metric on the 
plane? The geodesics of a Euclidean metric	
are the straight lines with respect to some coordinate system. Thus
our problem may be formulated as follows: when does there exist a coordinate
system $Y = Y(y)$ such that the solutions $y = y(t)$ of $(1)$ correspond
to straight lines $Y = Y(y(t)) = at+b$, $a,b \in {R}^{2}$? We therefore
seek to determine the existence of a transformation $\Psi:{R}^{2}\times
{R}\rightarrow  {R}^{2}\times {R}$ of the form
\[ \Psi (y,t) = (Y(y),t) \]
such that  $\Psi$ transforms the equations
\[  d^{2}Y^{1}/d^{2}t = 0  \hspace*{.3in} and \hspace*{.3in}
d^{2}Y^{2}/d^{2}t = 0 \hspace{1in} (2)  \]
into the equations $(1)$. Any transformation of the form $\Psi$ above transforms
$(2)$ into equations of the form
\[ d^{2}y^{i}/d^{2}t + \Gamma^{i}_{jk}(y)\dot{y}^{j}\dot{y}^{k} = 0. \]
Consequently, $(1)$ must necessarily be of this form.
The terms $\Gamma^{i}_{jk}(y)$ define a connection $\nabla$
on the surface and we see that the solutions to $(1)$ are locally the geodesics of a
Euclidean space if and only if $\nabla$ is flat.
The problem above  requires that the solutions to $(1)$ are already
parameterized in such a fashion  that only  a change in  the coordinates of the 
surface  is sufficient to straighten them out into lines. In this paper, we are
interested in whether the solutions to $(1)$ may be reparameterized so
as to be straight lines in some coordinate system. That is, do there exist
coordinates $Y = Y(y)$ and a reparameterization of time $T = T(y,t)$ such that
the solutions $y = y(t)$ of $(1)$ correspond to straight lines $Y = Y(T) = aT+b$, 
$a,b \in {R}^{2}$? We therefore 
seek to determine the existence of a transformation $\Phi:{R}^{2}\times
{R}\rightarrow  {R}^{2}\times {R}$ of the form
\[ \Phi (y,t) = (Y(y),T(y,t)) \]
such that  $\Phi$ transforms the equations 
\[  d^{2}Y^{1}/d^{2}T = 0  \hspace*{.3in} and \hspace*{.3in}
d^{2}Y^{2}/d^{2}T = 0 \]
into the equations $(1)$.

Conceivably, other reparameterization criteria could be considered as well. For instance,
one might investigate the more restricted transformation $T=T(t)$ where time is 
reparameterized in a manner independent of the point on the surface or the more general
$T=T(y, \dot{y},t)$. Here, we shall content ourselves with spacetime reparameterizations
only and defer the other cases to another time and place.

It is convenient to place the problem in a more general setting; we consider the
equivalence of two pairs of ordinary differential equations
\[ \frac{d^{2}y^{1}}{d^{2}t} = f(y, \dot{y},t) \hspace*{1in} 
\frac{d^{2}y^{2}}{d^{2}t} = g(y,\dot{y},t)  \hspace*{1in} (3) \]
and 
\[ \frac{d^{2}Y^{1}}{d^{2}T} = F(Y, \dot{Y},T) \hspace*{1in}
 \frac{d^{2}Y^{2}}{d^{2}T} = G(Y, \dot{Y},T). \hspace*{1in} (4) \]
under transformations of the form
\[ \Phi (y,t) = (Y(y),T(y,t)). \]
The case $F=G=0$ is solved in section 4. We obtain an e-structure on a 12-dimensional
space with constant torsion. The solutions to a pair of ODEs belonging to this 
equivalence class have symmetries given by the group of fractal-linear transformations
on the plane. In section 5 we make no restrictions on $F,G,f$ and $g$ and carry the 
equivalence through for the generic case.

A similar problem was studied by S.S. Chern [3]. He has considered the geometry of 
a system of second order ODEs
\[ \frac{d^{2}y^{i}}{d^{2}t} = f^{i}(y, \dot{y},t) 
\hspace*{1in} i=1,...,n ,\]
under transformations of the form
\[ \left\{ \begin{array}{lll}
Y & = & Y(y,t) \\
T & = & t. \end{array} \right. \]
Prior to Chern, the local behaviour of systems of second order ODEs
has been studied by M.D.D. Kosambi [6] and by E. Cartan [2].

\section{The Equivalence Problem Formulated}
The equations $(3)$
may be represented by the Pfaffian system 
\[ I = \left\{ \begin{array}{lll}
  dy^{1}-p^{1}dt & = & 0 \\
  dy^{2}-p^{2}dt & = & 0 \\
  dp^{1}-fdt & = & 0 \\
  dp^{2}-gdt & = & 0 \end{array} \right.  \]
on $U \subseteq {R}^{5}$. Similarly we represent the equations $(4)$  
by the Pfaffian system 

\[  J = \left\{ \begin{array}{lll}
  dY^{1}-P^{1}dT & = & 0 \\
  dY^{2}-P^{2}dT & = & 0 \\
  dP^{1}-FdT & = & 0 \\
  dP^{2}-GdT & = & 0 \end{array} \right. \] 

on $V\subseteq {R}^{5}$.

Form the coframes 
\[ \omega = \left\{ \begin{array}{lll}
\omega^{1} & = & dy^{1} \\
\omega^{2} & = & dy^{1}+pdy^{2} \\
\omega^{3} & = & dy^{2}+qdt \\
\omega^{4} & = & dp+hdt \\
\omega^{5} & = & dq+gdt \end{array} \right. \hspace*{.2in} and \hspace*{.2in}
\Omega = \left\{ \begin{array}{lll}
\Omega^{1} & = & dY^{1} \\
\Omega^{2} & = & dY^{1}+PdY^{2} \\
\Omega^{3} & = & dY^{2}+QdT \\
\Omega^{4} & = & dP+HdT \\
\Omega^{5} & = & dQ+GdT \end{array} \right. \]
where
\[ \left\{ \begin{array}{lll}
p & = & -p^{1}(p^{2})^{-1} \\
h & = & (p^{2})^{-2}(p^{2}f-p^{1}g) \\
q & = & -p^{2} \end{array} \right. \hspace*{.05in} and \hspace*{.05in}
\left\{ \begin{array}{lll}
P & = & -P^{1}(P^{2})^{-1} \\
H & = & (P^{2})^{-2}(P^{2}F-P^{1}G) \\
Q & = & -P^{2}. \end{array} \right. \]
Observe that $I$ (resp. $J$) is spanned by $\omega^{2},...,\omega^{5}$
(resp. $\Omega^{2},...,\Omega^{5}$).

Consider those transformations
\[ \Phi :U\rightarrow V \] 
such that
\[ (*) \hspace*{.5in} \left\{  \begin{array}{l}
(1) \hspace*{.1in} \Phi^{*}(J) = I , \hspace*{.05in}  and \\
(2) \hspace*{.1in} \Phi(y,p,t) = (Y(y), P(y,p,t),
      T(y,t)). \end{array} \right. \]
This is an {\em overdetermined} equivalence problem (cf. [5]). Nevertheless,
we shall not follow the approach given in the above reference in that all
the information contained in $(*)$ may be encoded by an apprpriately chosen
coframe and group, as follows.

Let $G$ be the subgroup of $GL(5,{R})$ whose elements are represented
by 
\[ \left(
\begin{tabular}{c|c|c}
$a$ & $\begin{array}{cc} b & 0 \end{array}$ & $0$ \\ \hline
$\begin{array}{c} 0 \\ 0 \end{array}$ & $\begin{array}{cc} c & 0 \\
e & f \end{array}$ & $0$ \\ \hline 
$0$ & $M$ & $N$ 
\end{tabular} 
\right). \]
It is easily shown that a diffeomerphism $\Phi :U\rightarrow V$ satisfies $(*)$
if and only if 
\[ \Phi^{*}(\Omega) = \gamma \omega \] 
for some  $\gamma:U\rightarrow G$. 
In order to avoid the awkward presence of the
map $\gamma$, form the lifted coframe $\eta$  (resp. $H$) on 
the  G-bundle $G\times U$ (resp. $G\times V$): 
\[ \eta := S\omega, \hspace*{.3in} (resp. \hspace*{.1in}  H := S\Omega), \]
where $S:G\rightarrow GL(5,{R})$ is the natural injection.
It can be shown that 
\[ \Phi^{*}(\Omega) = \gamma \omega \] 
if and only if 
\[ \bar{\Phi}^{*}(H) = \eta,\]
for some  $\bar{\Phi}:G\times U\rightarrow G\times V$. $\bar{\Phi}$
is related to $\Phi$ and $\gamma$ by
$\bar{\Phi}(g,u):= (g\gamma (u)^{-1},\Phi(u))$, for all $(g,u)\in G\times U$.
We arrive at the following formulation of the equivalence problem: \\
\\

{\it The two systems of ordinary differential equations (3) and (4) are equivalent
with respect to a diffeomorphism  $\Phi :U\rightarrow V$ of the form
$\Phi (y,t) = (Y(y),T(y,t))$  if and only if there exists a diffeomorphism
$\bar{\Phi}(g,u):G\times U\rightarrow G\times V$ satisfying
\[ \bar{\Phi}^{*}(H) = \eta.\] }
\\

With this characterization of equivalence we now proceed to reduce the group $G$.

\section{The First Normalization}
The structure equations, after making the obvious absorptions are 
\[ d\eta = \left( 
\begin{tabular}{c|c|c}
$\alpha$ & $\begin{array}{cc} \beta & 0 \end{array}$ & $0$ \\ \hline
$\begin{array}{c} 0 \\ 0 \end{array}$ & $\begin{array}{cc} \gamma & 0 \\
\epsilon & \phi  \end{array}$ & $0$ \\ \hline
$0$ & $\mu$ & $\nu$
\end{tabular} \right)\eta +
\left(
\begin{tabular}{c}
$0$ \\ \hline
$J\eta^{1}\eta^{3} + A_{1}\eta^{1}\eta^{4} + A_{2}\eta^{1}\eta^{5}$ \\
$B_{1}\eta^{1}\eta^{4}+B_{2}\eta^{1}\eta^{5}$ \\ \hline
$0$ \end{tabular} \right). \]

The equations 
\[ \begin{array} {lll}
0 & = & d^{2}\eta^{2}\eta^{2} \\
0 & = & d^{2}\eta^{3}\eta^{2}\eta^{3} 
\end{array} \hspace*{1in} (5)  \]
give us the following  infinitesimal group action on the torsion tensor.
\[ \left. \begin{array}{lll}
dJ + (\alpha-\gamma+\phi)J+A\mu_{2} &  \equiv &  0 \\
dA + (\alpha - \gamma)A + Av & \equiv & 0 \\
dB + (\alpha -\phi)B + Bv -\epsilon A & \equiv & 0 
\end{array} \right\} \hspace{.22in} mod \hspace{.1in} base.  \]
Therefore the group action on the torsion tensor is
\[ \left\{ \begin{array}{lll}
J & = & (J_{0}-A_{0}N^{-1}M_{2})a^{-1}cf^{-1} \\
A & = & a^{-1}cA_{0}N^{-1} \\
B & = & a^{-1}(fB_{0}+eA_{0})N^{-1}, \end{array} \right. \]
where $J_{0},A_{0},$ and $B_{0}$ denote the tensors $J,A$ and $B$, 
respectively, at the group identity. A parametric calculation shows that
\[ J_{0} = -(pq)^{-1}h, \hspace*{.2in} A_{0} = (p^{-1},0), \hspace*{.2in}
   B_{0}=(0,-(pq)^{-1}).\]
We may therefore normalize 
\[ J = 0, \hspace*{.1in}  A = (0,1), \hspace{.1in} and \hspace{.1in}
B=(1,0). \] 
Then
\[ \left. \begin{array}{lll}
\mu^{2}_{2} & \equiv & 0 \\
v^{2}_{1} & \equiv & 0 \\
\epsilon -v^{1}_{2} & \equiv & 0 \\
\alpha -\gamma + v^{2}_{2} & \equiv & 0 \\
\phi -\gamma  - v^{1}_{1} + v^{2}_{2} & \equiv & 0  \end{array} \right\}
\hspace{.22in} mod \hspace{.1in} base.  \]
Write
\[ \left\{ \begin{array}{lll}
\mu^{2}_{2} & = & A_{i}\eta^{i} \\
v^{2}_{1} & = & B_{i}\eta^{i} \\
\epsilon & = & v^{1}_{2} + C_{i}\eta^{i} \\
\alpha & = & \gamma -v^{2}_{2} + D_{i}\eta^{i} \\
\phi & = & \gamma + v^{1}_{1} -v^{2}_{2} +E_{i}\eta^{i} 
\end{array} \right. \]
Substituting these values  into (5) 
results in
\[ B_{3}=A_{4}, \hspace{.1in} D_{3}=A_{5}, \hspace{.1in} D_{4} = B_{5},
\hspace{.1in} E_{5} = C_{4}+D_{5}. \]
The structure equations become
\[ d\eta = \left(
\begin{tabular}{c|c|c}
$\gamma -v^{2}_{2}$ & $\begin{array}{cc} \beta \hspace*{.4in}  &
 0 \hspace*{.2in}  \end{array}$ & 
$\begin{array}{cc} 0 & 0 \end{array}$ \\ \hline
$\begin{array}{c} 0 \\ 0 \end{array}$ & $\begin{array}{cc} \hspace*{.1in} 
\gamma & 0 \\ \hspace*{.1in} v^{1}_{2} & \gamma +v^{1}_{1}-v^{2}_{2} 
 \end{array}$ & $\begin{array}{cc}
0 & 0 \\ 0 & 0 \end{array}$  \\ \hline
$\begin{array}{c} 0 \\ 0 \end{array}$ & $\begin{array}{cc} \mu^{1}_{1} \hspace*{.4in}
 & \mu^{1}_{2} \hspace*{.2in} \\ \mu^{2}_{1} \hspace*{.4in} 
 &  0 \hspace*{.2in} \end{array}$ & $\begin{array}{cc} v^{1}_{1} &
v^{1}_{2} \\ 0 & v^{2}_{2} \end{array}$
\end{tabular} \right)
\eta + \left(
\begin{tabular}{c}
$0$ \\ \hline $\eta^{1}\eta^{5}$ \\ $\eta^{1}\eta^{4}$ \\ \hline 
$\begin{array}{c} 0 \\ 0 \end{array}$
\end{tabular} \right) + T, \]
where 
\[ T = \left( 
\begin{tabular}{c}
$A_{5}\eta^{3}\eta^{1}+B_{5}\eta^{4}\eta^{1}+D_{5}\eta^{5}\eta^{1}$ \\ \hline
$0$ \\
$C_{1}\eta^{1}\eta^{2}+C_{3}\eta^{3}\eta^{2}+C_{4}\eta^{4}\eta^{2}+C_{5}\eta^{5}
\eta^{2}+E_{1}\eta^{1}\eta^{3}+E_{2}\eta^{2}\eta^{3}+E_{4}\eta^{4}\eta^{3}
+(C_{4}+D_{5})\eta^{5}\eta^{3}$ \\ \hline
$0$ \\
$A_{1}\eta^{1}\eta^{3}+A_{5}\eta^{5}\eta^{3}+B_{1}\eta^{1}
\eta^{4}+B_{5}\eta^{5}\eta^{4}$ 
\end{tabular} \right). \]
Absorb the torsion as follows:
\[ \begin{array}{lll}
v^{1}_{1} & \rightarrow & v^{1}_{1}+(E_{4}-B_{5})\eta^{4} + C_{4}\eta^{5} \\
v^{1}_{2} & \rightarrow & v^{1}_{2} + (C_{3}-E_{2})\eta^{3} + 
            C_{4}\eta^{4}+C_{5}\eta^{5} \\
v^{2}_{2} & \rightarrow & v^{2}_{2} -A_{5}\eta^{3}-B_{5}\eta^{4}-D_{5}\eta^{5} \\
\mu^{1}_{2} & \rightarrow & \mu^{1}_{2} + (C_{3}-E_{2})\eta^{5}.  \end{array} \]
The new structure equations are 
\[ d\eta = \left(
\begin{tabular}{c|c|c}
$\gamma -v^{2}_{2}$ & $\begin{array}{cc} \beta \hspace*{.4in}  &
 0 \hspace*{.2in}  \end{array}$ &
$\begin{array}{cc} 0 & 0 \end{array}$ \\ \hline
$\begin{array}{c} 0 \\ 0 \end{array}$ & $\begin{array}{cc} \hspace*{.1in}
\gamma & 0 \\ \hspace*{.1in} v^{1}_{2} & \gamma +v^{1}_{1}-v^{2}_{2}
 \end{array}$ & $\begin{array}{cc}
0 & 0 \\ 0 & 0 \end{array}$  \\ \hline
$\begin{array}{c} 0 \\ 0 \end{array}$ & $\begin{array}{cc} \mu^{1}_{1} \hspace*{.4in}
 & \mu^{1}_{2} \hspace*{.2in} \\ \mu^{2}_{1} \hspace*{.4in}
 &  0 \hspace*{.2in} \end{array}$ & $\begin{array}{cc} v^{1}_{1} &
v^{1}_{2} \\ 0 & v^{2}_{2} \end{array}$
\end{tabular} \right) \eta +
\left( \begin{tabular}{c}
$0$ \\ \hline
$\eta^{1}\eta^{5}$ \\
$\eta^{1}\eta^{4} + A\eta^{1}\eta^{2} + B\eta^{1}\eta^{3}$ \\ \hline
$0$ \\
$C\eta^{1}\eta^{3} + D\eta^{1}\eta^{4}$
\end{tabular} \right). \]

The equations
\[ \begin{array}{lll}
0 & = & d^{2}\eta^{5}\eta^{2}\eta^{5} \\
0 & = & d^{2}\eta^{3}\eta^{3}\eta^{4}\eta^{5} + d^{2}\eta^{4}\eta^{2}\eta^{3}
        \eta^{4} \\
0 & = & d^{2}\eta^{2}\eta^{3}\eta^{4}\eta^{5}+ d^{2}\eta^{3}\eta^{2}\eta^{4}\eta^{5}
+ d^{2}\eta^{4}\eta^{2}\eta^{3}\eta^{5}+ d^{2}\eta^{5}\eta^{2}\eta^{3}\eta^{4}
\end{array} \]
give the following
infinitesimal group action on the torsion tensor:
\[ \left. \begin{array}{lll}
dA+A(\gamma-v^{1}_{1})+ Bv^{1}_{2}+2\mu^{1}_{1} &  \equiv &  0 \\
dB + (\gamma-v^{2}_{2})B+2v^{1}_{2}D+2\mu^{1}_{2}-2\mu^{2}_{1} & \equiv & 0 \\
dC + (2\gamma + v^{1}_{1}-3v^{2}_{2})C + \mu^{1}_{2}D & \equiv & 0 \\
dD+(\gamma + v^{1}_{1} - 2v^{2}_{2})D & \equiv & 0
\end{array} \right\} \hspace{.22in} mod \hspace{.1in} base.  \]

\section{Geodesics of Flat, Symmetric Connections}
In this section the equivalence problem is carried through for the
case $F=G=0$. This will lead to an e-structure with constant torsion on 
a 12-dimensional space. The only invariants are therefore constant 
invariants and thus we obtain a complete solution to the problem of 
equivalence.  

\subsection{The Second Normalization}	
A parametric calculation will show that at the identity $A_{0}=B_{0}=C_{0}=D_{0}=0$,
hence $D\equiv C\equiv 0$. This leaves the following two equations:
\[ \left. \begin{array}{lll}
dA+A(\gamma-v^{1}_{1})+Bv^{1}_{2}+2\mu^{1}_{1} &  \equiv &  0 \\
dB + (\gamma-v^{2}_{2})B+2\mu^{1}_{2}-2\mu^{2}_{1} & \equiv & 0
\end{array} \right\} \hspace{.22in} mod \hspace{.1in} base.  \]
Normalize $A\equiv  B\equiv 0$. We then have
\[ \left.  \begin{array}{lll}
\mu^{1}_{1} & \equiv & 0 \\
\mu^{1}_{2} & \equiv & \mu^{2}_{1} \end{array} \right\}
\hspace{.22in} mod \hspace{.1in} base.  \]
This produces new torsion by
\[ \begin{array}{lll}
\mu^{1}_{1} & = & A_{i}\eta^{i} \\
\mu^{1}_{2} & = & \mu^{2}_{1}+B_{i}\eta^{i}. \end{array} \]
This gives the following structure equations:
\[ d\eta = \left(
\begin{tabular}{c|c|c}
$\gamma -v^{2}_{2}$ & $\begin{array}{cc} \beta \hspace*{.4in}  &
 0 \hspace*{.2in}  \end{array}$ &
$\begin{array}{cc} 0 & 0 \end{array}$ \\ \hline
$\begin{array}{c} 0 \\ 0 \end{array}$ & $\begin{array}{cc} \hspace*{.1in}
\gamma & 0 \\ \hspace*{.1in} v^{1}_{2} & \gamma +v^{1}_{1}-v^{2}_{2}
 \end{array}$ & $\begin{array}{cc}
0 & 0 \\ 0 & 0 \end{array}$  \\ \hline
$\begin{array}{c} 0 \\ 0 \end{array}$ & $\begin{array}{cc} 0 \hspace*{.4in}
 & \mu^{2}_{1} \hspace*{.2in} \\ \mu^{2}_{1} \hspace*{.4in}
 &  0 \hspace*{.2in} \end{array}$ & $\begin{array}{cc} v^{1}_{1} &
v^{1}_{2} \\ 0 & v^{2}_{2} \end{array}$
\end{tabular} \right) \eta +
\left( \begin{tabular}{c}
$0$ \\ \hline
$\eta^{1}\eta^{5}$ \\
$\eta^{1}\eta^{4} $ \\ \hline
$0$ \\
$0$
\end{tabular} \right) +T, \]
where
\[ T = \left( \begin{array}{c}
0 \\
0 \\
0 \\
A_{1}\eta^{1}\eta^{2}+A_{3}\eta^{3}\eta^{2}+A_{4}\eta^{4}\eta^{2}+
A_{5}\eta^{5}\eta^{2}+ B_{1}\eta^{1}\eta^{3}+B_{2}\eta^{2}\eta^{3}+
B_{4}\eta^{4}\eta^{3}+B_{5}\eta^{5}\eta^{3} \\
0 
\end{array} \right). \]

Absorb
\[ \begin{array}{lll}
v^{1}_{1} & \rightarrow & v^{1}_{1}-A_{4}\eta^{2}-B_{4}\eta^{3} \\
v^{1}_{2} & \rightarrow & v^{1}_{2} -A_{5}\eta^{2} -B_{5}\eta^{3} \\
\mu^{2}_{1} & \rightarrow & \mu^{2}_{1}-A_{3}\eta^{2}+B_{2}\eta^{2} \\
\gamma & \rightarrow & \gamma + A_{4}\eta^{2}-B_{5}\eta^{2} \\
\beta & \rightarrow & \beta + A_{4}\eta^{1}-B_{5}\eta^{1}. \end{array} \]
We then obtain
\[ d\eta = \left(
\begin{tabular}{c|c|c}
$\gamma -v^{2}_{2}$ & $\begin{array}{cc} \beta \hspace*{.4in}  &
 0 \hspace*{.2in}  \end{array}$ &
$\begin{array}{cc} 0 & 0 \end{array}$ \\ \hline
$\begin{array}{c} 0 \\ 0 \end{array}$ & $\begin{array}{cc} \hspace*{.1in}
\gamma & 0 \\ \hspace*{.1in} v^{1}_{2} & \gamma +v^{1}_{1}-v^{2}_{2}
 \end{array}$ & $\begin{array}{cc}
0 & 0 \\ 0 & 0 \end{array}$  \\ \hline
$\begin{array}{c} 0 \\ 0 \end{array}$ & $\begin{array}{cc} 0 \hspace*{.4in}
 & \mu^{2}_{1} \hspace*{.2in} \\ \mu^{2}_{1} \hspace*{.4in}
 &  0 \hspace*{.2in} \end{array}$ & $\begin{array}{cc} v^{1}_{1} &
v^{1}_{2} \\ 0 & v^{2}_{2} \end{array}$
\end{tabular} \right) \eta +
\left( \begin{tabular}{c}
$0$ \\ \hline
$\eta^{1}\eta^{5}$ \\
$\eta^{1}\eta^{4} $ \\ \hline
$A\eta^{1}\eta^{2}+B\eta^{1}\eta^{3}$ \\
$0$
\end{tabular} \right).  \]
The equations
\[ \begin{array}{lll}
0 & = & d^{2}\eta^{4}\eta^{3}\eta^{4}\eta^{5} \\
0 & = & d^{2}\eta^{4}\eta^{2}\eta^{4}\eta^{5}
       + d^{2}\eta^{5}\eta^{3}\eta^{4}\eta^{5}  \end{array} \]
give
\[ \left. \begin{array}{lll}
dA + (2\gamma-v^{1}_{1}-v^{2}_{2})A + Bv^{1}_{2} & \equiv & 0 \\
dB+(2\gamma-2v^{2}_{2})B & \equiv & 0 \end{array} \right\}
\hspace{.22in} mod \hspace{.1in} base.  \]
Now at the identity $A_{0}=B_{0}=0$. Thus $A=B=0$, and hence
\[ d\eta = \left(
\begin{tabular}{c|c|c}
$\gamma -v^{2}_{2}$ & $\begin{array}{cc} \beta \hspace*{.4in}  &
 0 \hspace*{.2in}  \end{array}$ &
$\begin{array}{cc} 0 & 0 \end{array}$ \\ \hline
$\begin{array}{c} 0 \\ 0 \end{array}$ & $\begin{array}{cc} \hspace*{.1in}
\gamma & 0 \\ \hspace*{.1in} v^{1}_{2} & \gamma +v^{1}_{1}-v^{2}_{2}
 \end{array}$ & $\begin{array}{cc}
0 & 0 \\ 0 & 0 \end{array}$  \\ \hline
$\begin{array}{c} 0 \\ 0 \end{array}$ & $\begin{array}{cc} 0 \hspace*{.4in}
 & \mu^{2}_{1} \hspace*{.2in} \\ \mu^{2}_{1} \hspace*{.4in}
 &  0 \hspace*{.2in} \end{array}$ & $\begin{array}{cc} v^{1}_{1} &
v^{1}_{2} \\ 0 & v^{2}_{2} \end{array}$
\end{tabular} \right) \eta +
\left( \begin{tabular}{c}
$0$ \\ \hline
$\eta^{1}\eta^{5}$ \\
$\eta^{1}\eta^{4} $ \\ \hline
$0$ \\
$0$
\end{tabular} \right).  \]
Only constant torsion remains so the system must be prolonged.

\subsection{Prolongation}
Now dim ${\cal {G}}^{(1)} = 2$.
The first prolongation corresponds to the following arbitrariness in the
tableau:
\[  \begin{array}{lll}
\theta^{1} & := & v^{1}_{1} \\
\theta^{2} & := & v^{1}_{2}+a\eta^{3} \\
\theta^{3} & := & v^{2}_{2}+ a\eta^{2} \\
\theta^{4} & := & \gamma +2a\eta^{2} \\
\theta^{5} & := & \beta +a\eta^{1}+b\eta^{2} \\
\theta^{6} & := & \mu^{2}_{1} +a\eta^{5} \end{array} \]
where $a,b \in  {R}$ are arbitrary.

The structure equations for $\eta$ may be written,
\[ d\eta = \left(
\begin{tabular}{c|c|c}
$\theta^{4}-\theta^{3}$ & $\begin{array}{cc} \theta^{5} \hspace*{.4in}  &
 0 \hspace*{.2in}  \end{array}$ &
$\begin{array}{cc} 0 & 0 \end{array}$ \\ \hline
$\begin{array}{c} 0 \\ 0 \end{array}$ & $\begin{array}{cc} \hspace*{.1in}
\theta^{4} & 0 \\ \hspace*{.1in} \theta^{2} & \theta^{4}+\theta^{1}-\theta^{3}
 \end{array}$ & $\begin{array}{cc}
0 & 0 \\ 0 & 0 \end{array}$  \\ \hline
$\begin{array}{c} 0 \\ 0 \end{array}$ & $\begin{array}{cc} 0 \hspace*{.4in}
 & \theta^{6} \hspace*{.2in} \\ \theta^{6} \hspace*{.4in}
 &  0 \hspace*{.2in} \end{array}$ & $\begin{array}{cc} \theta^{1} &
\theta^{1} \\ 0 & \theta^{3} \end{array}$
\end{tabular} \right) \eta +
\left( \begin{tabular}{c}
$0$ \\ \hline
$\eta^{1}\eta^{5}$ \\
$\eta^{1}\eta^{4} $ \\ \hline
$0$ \\
$0$
\end{tabular} \right).  \]

Taking $d^{2}$ of the above structure equations we obtain
the following:
\[ \begin{array}{lll}
0 & = & -d\theta^{3}\eta^{1}+d\theta^{4}\eta^{1}+d\theta^{5}\eta^{2}+
          \theta^{3}\theta^{5}\eta^{2}-\theta^{5}\eta^{1}\eta^{5} \\
0 & = & d\theta^{4}\eta^{2}+\theta^{5}\eta^{2}\eta^{5}+\theta^{6}\eta^{1}\eta^{2} \\
0 & = & d\theta^{1}\eta^{3}+d\theta^{2}\eta^{2}-d\theta^{3}\eta^{3}+
        d\theta^{4}\eta^{3}-\theta^{1}\theta^{2}\eta^{2}+\theta^{3}\theta^{2}\eta^{2}
        +\theta^{5}\eta^{2}\eta^{4}+\theta^{6}\eta^{1}\eta^{3} \\
0 & = & d\theta^{1}\eta^{4}+d\theta^{2}\eta^{5}+d\theta^{6}\eta^{3}
        -\theta^{2}\theta^{3}\eta^{5}+\theta^{6}\theta^{3}\eta^{3}+
        \theta^{4}\theta^{6}\eta^{3}-\theta^{6}\eta^{1}\eta^{4}-
         \theta^{1}\theta^{2}\eta^{5} \\
0 & = & d\theta^{3}\eta^{5}+d\theta^{6}\eta^{2}-\theta^{3}\theta^{6}\eta^{2}
        +\theta^{4}\theta^{6}\eta^{2}-\theta^{6}\eta^{1}\eta^{5}.
\end{array} \]
It follows from a somewhat lengthy but straightforward calculation that
the structure equations for $\theta$ are given by 
\[ \begin{array}{lll}
d\theta^{1} & = & A_{1}\eta^{5}\eta^{2}+A_{2}\eta^{3}\eta^{2}+A_{3}\eta^{4}\eta^{2}+
         A_{4}\eta^{4}\eta^{3}+A_{5}\eta^{5}\eta^{3}+\theta^{5}\eta^{5}+
         \theta^{6}\eta^{1} \\
d\theta^{2} & = & \Phi^{1}\eta^{3} + A_{1}\eta^{4}\eta^{2}+B\eta^{5}\eta^{2}+
         A_{5}\eta^{4}\eta^{3}+\theta^{1}\theta^{2}-\theta^{3}\theta^{2}+
         \theta^{5}\eta^{4} \\
d\theta^{3} & = & \Phi^{1}\eta^{2} + C\eta^{5}\eta^{2}+2\theta^{5}\eta^{5}+
         \theta^{6}\eta^{1} \\
d\theta^{4} & = & 2\Phi^{1}\eta^{2}+D_{1}\eta^{3}\eta^{2}+D_{2}\eta^{4}\eta^{2}
         +(C-A_{1})\eta^{5}\eta^{2}+\theta^{5}\eta^{5}-\theta^{6}\eta^{1} \\
d\theta^{5} & = & \Phi^{1}\eta^{1}+\Phi^{2}\eta^{2}+D_{1}\eta^{3}\eta^{1}+
          D_{2}\eta^{4}\eta^{1}-A_{1}\eta^{5}\eta^{1}-\theta^{3}\theta^{5} \\
d\theta^{6} & = & \Phi^{1}\eta^{5}+F\eta^{3}\eta^{2}+A_{2}\eta^{4}\eta^{2}
           +\theta^{3}\theta^{6}-\theta^{4}\theta^{6}, \end{array} \]
where $\Phi^{1}$ and $\Phi^{2}$ are 1-forms in the new tableau.
The equations
\[ \begin{array}{lll}
0 & = & d^{2}\theta^{1} \\
0 & = & d^{2}\theta^{2}\eta^{3} \\
0 & = & d^{2}\theta^{2}\eta^{2}+d^{2}\theta^{3}\eta^{3} \\
0 & = & 2d^{2}\theta^{3}-d^{2}\theta^{4} \\
0 & = & d^{2}\theta^{6}\eta^{5} \end{array} \]
give the following infinitesimal action on the torsion tensor
\[ \left. \begin{array}{lll}
dA_{1} -\Phi^{2} & \equiv & 0 \\
dA_{2} & \equiv & 0 \\
dA_{3} & \equiv & 0 \\
dA_{4} & \equiv & 0 \\
dA_{5} & \equiv & 0 \\
dB & \equiv & 0 \\
dC -2\Phi ^{2} & \equiv & 0 \\
dD_{1} & \equiv & 0 \\
dD_{2} & \equiv & 0 \\
dF & \equiv & 0 \end{array} \right\}
\hspace{.22in} mod \hspace{.1in} base.  \]
At the group identity,
\[ A_{2}=A_{3}=A_{4}=A_{5}=B=D_{1}=D_{2}=F=0, \]
therefore these terms are identically zero.
Thus,
\[ (dA_{1}+A_{1}(\theta^{3}+\theta^{4})-\Phi^{2})\eta^{4}\eta^{2}\eta^{3}
  + (2A_{1}-C)\theta^{2}\eta^{5}\eta^{2}\eta^{3} =0. \]
Consequently, \[ C=2A_{1}. \]
We also have,
\[ dA_{1}+A_{1}(\theta^{3}+\theta^{4})-\Phi^{2}  \equiv 0
\hspace{.22in} mod \hspace{.1in} \eta^{2},\eta^{5}.  \]
Therefore,
\[ dA_{1}+A_{1}(\theta^{3}+\theta^{4})-\Phi^{2} + l\eta^{2} =0 \]
for some function $l$.
We obtain the following structure equations
\[ \begin{array}{lll}
d\theta^{1} & = & A\eta^{5}\eta^{2} +\theta^{5}\eta^{5}+ \theta^{6}\eta^{1} \\
d\theta^{2} & = & \Phi^{1}\eta^{3} + A\eta^{4}\eta^{2}
         +\theta^{1}\theta^{2}-\theta^{3}\theta^{2}+ \theta^{5}\eta^{4} \\
d\theta^{3} & = & \Phi^{1}\eta^{2} + 2A\eta^{5}\eta^{2}+2\theta^{5}\eta^{5}+
         \theta^{6}\eta^{1} \\
d\theta^{4} & = & 2\Phi^{1}\eta^{2}
         +A\eta^{5}\eta^{2}+\theta^{5}\eta^{5}-\theta^{6}\eta^{1} \\
d\theta^{5} & = & \Phi^{1}\eta^{1}+\Phi^{2}\eta^{2}
          -A\eta^{5}\eta^{1}-\theta^{3}\theta^{5} \\
d\theta^{6} & = & \Phi^{1}\eta^{5}
           +\theta^{3}\theta^{6}-\theta^{4}\theta^{6}, \end{array} \]
where $A$ has been written for $A_{1}$.

\subsection{The Third Normalization}
We have the infinitesimal group action on the torsion given by
\[ dA_{1}+A_{1}(\theta^{3}+\theta^{4})-\Phi^{2}  \equiv 0
\hspace{.22in} mod \hspace{.1in} base.  \]
Normalize $A=0$. Then $\Phi^{2}=l\eta^{2}$ and hence
\[ \begin{array}{lll}
d\theta^{1} & = & \theta^{5}\eta^{5}+ \theta^{6}\eta^{1} \\
d\theta^{2} & = & \Phi^{1}\eta^{3} +
         \theta^{1}\theta^{2}-\theta^{3}\theta^{2}+ \theta^{5}\eta^{4} \\
d\theta^{3} & = & \Phi^{1}\eta^{2} +2\theta^{5}\eta^{5}+
         \theta^{6}\eta^{1} \\
d\theta^{4} & = & 2\Phi^{1}\eta^{2}
         +\theta^{5}\eta^{5}-\theta^{6}\eta^{1} \\
d\theta^{5} & = & \Phi^{1}\eta^{1}
          -\theta^{3}\theta^{5} \\
d\theta^{6} & = & \Phi^{1}\eta^{5}
           +\theta^{3}\theta^{6}-\theta^{4}\theta^{6}. \end{array} \]
We have constant torsion and an e-structure. Thus
$(\eta^{1},...,\eta^{5},\theta^{1},...,\theta^{6},\Phi)$ is an invariant
coframe (here we have written $\Phi = \Phi^{1}$.)
The equations
\[ (d\Phi -\Phi\theta^{4}-\theta^{5}\theta^{6})\eta^{i} = 0  \]
for $i=1,2,3$ and $5$ give
\[ d\Phi = \Phi\theta^{4}+\theta^{5}\eta^{6}. \]
Thus the structure equations are
\[ \left. \begin{array}{lll}
d\eta^{1} & = & (\theta^{4}-\theta^{3})\eta^{1}+\theta^{5}\eta^{2} \\
d\eta^{2} & = & \theta^{4}\eta^{2}+\eta^{1}\eta^{5} \\
d\eta^{3} & = & \theta^{2}\eta^{2}+(\theta^{4}+\theta^{1}-\theta^{3})\eta^{3}
                +\eta^{1}\eta^{4} \\
d\eta^{4} & = & \theta^{6}\eta^{3}+\theta^{1}\eta^{4}+\theta^{2}\eta^{5} \\
d\eta^{5} & = & \theta^{6}\eta^{2}+\theta^{3}\eta^{5} \\
d\theta^{1} & = & \theta^{5}\eta^{5}+\theta^{6}\eta^{1} \\
d\theta^{2} & = & \Phi\eta^{3}+(\theta^{1}-\theta^{3})\theta^{2}+
                  \theta^{5}\eta^{4} \\
d\theta^{3} & = & \Phi\eta^{2}+2\theta^{5}\eta^{5}+\theta^{6}\eta^{1} \\
d\theta^{4} & = & 2\Phi\eta^{2}+\theta^{5}\eta^{5}-\theta^{6}\eta^{1} \\
d\theta^{5} & = & \Phi\eta^{1}-\theta^{3}\theta^{5} \\
d\theta^{6} & = & \Phi\eta^{5}+(\theta^{3}-\theta^{4})\theta^{6} \\
d\Phi & = & \Phi\theta^{4}+\theta^{5}\theta^{6}.
\end{array} \right\} \hspace*{.5in} (6)  \]

We have shown the following:
\begin{theorem}
The solutions $y(t)$ of 
\[ d^{2}y^{1}/d^{2}t = f(y,\dot{y},t)  \hspace*{.3in} and \hspace*{.3in}
d^{2}y^{2}/d^{2}t = g(y,\dot{y},t),  \hspace{1in} (7) \]
are the geodesics of a Euclidean space
with respect to a transformation of the form $(Y,T)=(Y(y),T(y,t))$
if and only if equations $(7)$ yield the structure equations $(6)$.
\end{theorem}

The structure equations $(6)$ are the Mauer-Cartan equations for a Lie-group:
the group of fractal-linear transformations ${\cal {F}}$ on the plane. Recall
that the fractal-linear transformations on the plane are those transformations
$A$ of the form
\[ A = (\bar{y},\bar{t}) = (\bar{y}(y),\bar{t}(y,t)) \]
where
\[ \bar{y}^{1}(y)  =  \frac{b^{1}_{0}+b^{1}_{1}y^{1}+b^{1}_{2}y^{2}}
                         {a_{0}+a_{1}y^{1}+a_{2}y^{2}} \]
\[ \bar{y}^{2}(y)  =  \frac{b^{2}_{0}+b^{2}_{1}y^{1}+b^{2}_{2}y^{2}}
                         {a_{0}+a_{1}y^{1}+a_{2}y^{2}} \]
\[ \hspace*{.2in} \bar{t}(y,t)  =  \frac{t + c_{0}+c_{1}y^{1}+c_{2}y^{2}}
                         {a_{0}+a_{1}y^{1}+a_{2}y^{2}}, \]
where $a_{i},a^{i}_{j},b^{i}_{j}\in {R}$ are constants.
Let ${\cal {G}}$ denote the subgroup of $GL(4,{R})$  consisting 
of those invertible matrices $M$ whose first column is $^{t}(1,0,0,0)$. 
We may identify the transformation $A$ with the element
$\phi(A)$ in  ${\cal {G}}$ by
\[ \phi(A):= \left( \begin{tabular}{c|c|c}
$1$ & $c_{0}$ & $\begin{array}{cc} c_{1} & c_{2} \end{array}$ \\ \hline
$0$ & $a_{0}$ & $\begin{array}{cc} a_{1} & a_{2} \end{array}$ \\ \hline
$\begin{array}{c} 0 \\ 0 \end{array}$ & $\begin{array}{c}
b^{1}_{0} \\ b^{2}_{0} \end{array}$ & $\begin{array}{cc}
b^{1}_{1} & b^{1}_{2} \\ b^{2}_{1} & b^{2}_{2} \end{array}$ 
\end{tabular} \right). \] 
The map $\phi:{\cal {F}}\rightarrow {\cal {G}}$ is a group isomorphism
and so we may identify the fractal-linear transformations with ${\cal {G}}$.

\section{The Generic Case}
In this section we continue the equivalence problem from the end of
the first normalization (section 3).

\subsection{The Second Normalization}
Recall that after the first normalization we have the structure equations
\[ d\eta = \left(
\begin{tabular}{c|c|c}
$\gamma -v^{2}_{2}$ & $\begin{array}{cc} \beta \hspace*{.4in}  &
 0 \hspace*{.2in}  \end{array}$ &
$\begin{array}{cc} 0 & 0 \end{array}$ \\ \hline
$\begin{array}{c} 0 \\ 0 \end{array}$ & $\begin{array}{cc} \hspace*{.1in}
\gamma & 0 \\ \hspace*{.1in} v^{1}_{2} & \gamma +v^{1}_{1}-v^{2}_{2}
 \end{array}$ & $\begin{array}{cc}
0 & 0 \\ 0 & 0 \end{array}$  \\ \hline
$\begin{array}{c} 0 \\ 0 \end{array}$ & $\begin{array}{cc} \mu^{1}_{1} \hspace*{.4in}
 & \mu^{1}_{2} \hspace*{.2in} \\ \mu^{2}_{1} \hspace*{.4in}
 &  0 \hspace*{.2in} \end{array}$ & $\begin{array}{cc} v^{1}_{1} &
v^{1}_{2} \\ 0 & v^{2}_{2} \end{array}$
\end{tabular} \right) \eta +
\left( \begin{tabular}{c}
$0$ \\ \hline
$\eta^{1}\eta^{5}$ \\
$\eta^{1}\eta^{4} + A\eta^{1}\eta^{2} + B\eta^{1}\eta^{3}$ \\ \hline
$0$ \\
$C\eta^{1}\eta^{3} + D\eta^{1}\eta^{4}$
\end{tabular} \right) \]
with  infinitesimal group action on the torsion tensor given by
\[ \left. \begin{array}{lll}
dA+A(\gamma-v^{1}_{1})+ Bv^{1}_{2}+2\mu^{1}_{1} &  \equiv &  0 \\
dB + (\gamma-v^{2}_{2})B+2v^{1}_{2}D+2\mu^{1}_{2}-2\mu^{2}_{1} & \equiv & 0 \\
dC + (2\gamma + v^{1}_{1}-3v^{2}_{2})C + \mu^{1}_{2}D & \equiv & 0 \\
dD+(\gamma + v^{1}_{1} - 2v^{2}_{2})D & \equiv & 0
\end{array} \right\} \hspace{.22in} mod \hspace{.1in} base.  \]
In the generic case, $D\neq 0$.  
We thus normalize
\[ A=0, \hspace{.1in} B=0, \hspace{.1in} C=0, \hspace{.1in} D=1.\]
This gives us 
\[ \left. \begin{array}{lll}
\mu^{1}_{1} & \equiv & 0 \\
\mu^{1}_{2} & \equiv & 0 \\
\mu^{2}_{1} & \equiv & v^{1}_{2} \\
\gamma  & \equiv & 2v^{2}_{2}-v^{1}_{1}  \end{array} \right\}
\hspace{.22in} mod \hspace{.1in} base.  \]
Let 
\[  \begin{array}{lll}
\mu^{1}_{1} & =& A_{i}\eta^{i} \\
\mu^{1}_{2} & = & B_{i}\eta^{i} \\
\mu^{2}_{1} & = & v^{1}_{2} + C_{i}\eta^{i}\\
\gamma  & = & 2v^{2}_{2}-v^{1}_{1} + D_{i}\eta^{i}. \end{array} \]
It follows that \[ B_{4}=D_{3}. \]
The new structure equations are:
\[ d\eta = \left(
\begin{tabular}{c|c|c}
$v^{2}_{2}-v^{1}_{1}$ & $\begin{array}{cc}  \hspace*{.2in} \beta \hspace*{.2in}  &
\hspace*{.1in} 0 \hspace*{.2in}  \end{array}$ &
$\begin{array}{cc} 0 & 0 \end{array}$ \\ \hline
$\begin{array}{c} 0 \\ 0 \end{array}$ & $\begin{array}{cc} 
2v^{2}_{2}-v^{1}_{1} & 0 \hspace*{.2in} \\ v^{1}_{2} & v^{2}_{2} \hspace*{.2in}
 \end{array}$ & $\begin{array}{cc}
0 & 0 \\ 0 & 0 \end{array}$  \\ \hline
$\begin{array}{c} 0 \\ 0 \end{array}$ & $\begin{array}{cc} 
\hspace*{.2in} 0 \hspace*{.2in}
\hspace*{.1in} & 0 \hspace*{.2in} \\ \hspace*{.2in} v^{1}_{2} \hspace*{.2in}
 & \hspace*{.05in} 0 \hspace*{.2in} \end{array}$ & $\begin{array}{cc} v^{1}_{1} &
v^{1}_{2} \\ 0 & v^{2}_{2} \end{array}$
\end{tabular} \right)
\eta + \left(
\begin{tabular}{c}
$0$ \\ \hline $\eta^{1}\eta^{5}$ \\ $\eta^{1}\eta^{4}$ \\ \hline
$\begin{array}{c} 0 \\ \eta^{1}\eta^{4}  \end{array}$
\end{tabular} \right) + T, \]
where 
\[ T = \left( \begin{tabular}{c}
$B_{4}\eta^{3}\eta^{1}+D_{4}\eta^{4}\eta^{1}+D_{5}\eta^{5}\eta^{1}+D_{2}\eta^{2}
\eta^{1}$ \\ \hline
$D_{1}\eta^{1}\eta^{2}+B_{4}\eta^{3}\eta^{2}+D_{4}\eta^{4}\eta^{2}+D_{5}\eta^{5}
\eta^{2}$ \\
$D_{1}\eta^{1}\eta^{3}+D_{2}\eta^{2}\eta^{3}+D_{4}\eta^{4}\eta^{3}+D_{5}\eta^{5}
\eta^{3}$ \\ \hline
$A_{1}\eta^{1}\eta^{2}+A_{3}\eta^{3}\eta^{2}+A_{4}\eta^{4}\eta^{2}+A_{5}\eta^{5}
\eta^{2}+
B_{1}\eta^{1}\eta^{3}+B_{2}\eta^{2}\eta^{3}+B_{4}\eta^{4}\eta^{3}+B_{5}\eta^{5}
\eta^{3}$ \\
$C_{1}\eta^{1}\eta^{2}+C_{3}\eta^{3}\eta^{2}+C_{4}\eta^{4}\eta^{2}+C_{5}\eta^{5}
\eta^{2}$ \end{tabular} \right). \]
Absorb the torsion:
\[ \begin{array}{lll}
v^{1}_{1} & \rightarrow & v^{1}_{1}-B_{4}\eta^{3}-D_{4}\eta^{4}-A_{4}\eta^{2} \\
v^{1}_{2} & \rightarrow & v^{1}_{2} -A_{5}\eta^{2}+C_{3}\eta^{3} \\
\beta & \rightarrow & \beta + A_{4}\eta^{1}+C_{3}\eta^{1} \\
v^{2}_{2} & \rightarrow & v^{2}_{2} + (C_{3}+D_{2})\eta^{2}. \end{array} \]
The structure equations are then
\[ d\eta = \left(
\begin{tabular}{c|c|c}
$v^{2}_{2}-v^{1}_{1}$ & $\begin{array}{cc}  \hspace*{.2in} \beta \hspace*{.2in}  &
\hspace*{.1in} 0 \hspace*{.2in}  \end{array}$ &
$\begin{array}{cc} 0 & 0 \end{array}$ \\ \hline
$\begin{array}{c} 0 \\ 0 \end{array}$ & $\begin{array}{cc}
2v^{2}_{2}-v^{1}_{1} & 0 \hspace*{.2in} \\ v^{1}_{2} & v^{2}_{2} \hspace*{.2in}
 \end{array}$ & $\begin{array}{cc}
0 & 0 \\ 0 & 0 \end{array}$  \\ \hline
$\begin{array}{c} 0 \\ 0 \end{array}$ & $\begin{array}{cc}
\hspace*{.2in} 0 \hspace*{.2in}
\hspace*{.1in} & 0 \hspace*{.2in} \\ \hspace*{.2in} v^{1}_{2} \hspace*{.2in}
 & \hspace*{.05in} 0 \hspace*{.2in} \end{array}$ & $\begin{array}{cc} v^{1}_{1} &
v^{1}_{2} \\ 0 & v^{2}_{2} \end{array}$
\end{tabular} \right)
\eta + \left(
\begin{tabular}{c}
$A\eta^{5}\eta^{1}$ \\ \hline 
$\eta^{1}\eta^{5}+B\eta^{1}\eta^{2}+A\eta^{5}\eta^{2}$ \\
$\eta^{1}\eta^{4}+B\eta^{1}\eta^{3}+C\eta^{4}\eta^{3}+A\eta^{5}\eta^{3}$ \\ \hline
$D\eta^{1}\eta^{2}+E\eta^{3}\eta^{2}+F\eta^{1}\eta^{3}+G\eta^{5}\eta^{3}$ \\
$\eta^{1}\eta^{4}+H\eta^{1}\eta^{2}+I\eta^{4}\eta^{2}+J\eta^{5}\eta^{2}$
\end{tabular} \right). \]

\subsection{The Third Normalization}
The equations
\[  \begin{array}{lll}
0 & = & d^{2}\eta^{1}\eta^{2}\eta^{3}+d^{2}\eta^{2}\eta^{1}\eta^{3}+
        d^{2}\eta^{3}\eta^{1}\eta^{2} \\
0 & = & d^{2}\eta^{3}\eta^{2}\eta^{5}+d^{2}\eta^{5}\eta^{2}\eta^{3} \\
0 & = & d^{2}\eta^{4}\eta^{4}\eta^{5} \\
0 & = & d^{2}\eta^{4}\eta^{2}\eta^{4}-d^{2}\eta^{5}\eta^{4}\eta^{5} \\
0 & = & d^{2}\eta^{3}\eta^{3}\eta^{5}-d^{2}\eta^{5}\eta^{3}\eta^{5} \\
0 & = & -d^{2}\eta^{3}\eta^{4}\eta^{5}-d^{2}\eta^{3}\eta^{3}\eta^{4} +
        d^{2}\eta^{4}\eta^{2}\eta^{4}+d^{2}\eta^{5}\eta^{3}\eta^{4}
\end{array} \hspace*{1in} (8) \]
give us the following infinitesimal group action on the torsion tensor:
\[ \left. \begin{array}{lll}
dA+Av^{2}_{2}+Cv^{1}_{2}+2\beta & \equiv & 0 \\
dB+B(v^{2}_{2}-v^{1}_{1}) & \equiv & 0 \\
dC+Cv^{1}_{1} & \equiv & 0 \\
dD+3D(v^{2}_{2}-v^{1}_{1})+Fv^{1}_{2}-Hv^{1}_{2} & \equiv & 0 \\
dE+E(3v^{2}_{2}-2v^{1}_{1})-F\beta-Gv^{1}_{2} & \equiv & 0 \\
dF+2F(v^{2}_{2}-v^{1}_{1}) & \equiv & 0 \\
dG+G(2v^{2}_{2}-v^{1}_{1}) & \equiv & 0 \\
dH+2H(v^{2}_{2}-v^{1}_{1})-Bv^{1}_{2} & \equiv & 0 \\
dI+Iv^{2}_{2}-Cv^{1}_{2} & \equiv & 0 \\
dJ+J(2v^{2}_{2}-v^{1}_{1})+B\beta + Iv^{1}_{2} & \equiv & 0 \end{array}
\right\} 
\hspace{.22in} mod \hspace{.1in} base.  \]
At this point, the group is the subgroup of $GL(5,{R})$ whose elenents
are of the form
\[ \left(
\begin{tabular}{c|c|c}
$ca^{-1}$ & $\begin{array}{cc}  e &  0  \end{array}$ &
$\begin{array}{cc} 0 & 0 \end{array}$ \\ \hline
$\begin{array}{c} 0 \\ 0 \end{array}$ & $\begin{array}{cc}
c^{2}a^{-1} & 0 \hspace*{.2in} \\ a^{-1}bc & c
 \end{array}$ & $\begin{array}{cc}
0 & 0 \\ 0 & 0 \end{array}$  \\ \hline
$\begin{array}{c} 0 \\ 0 \end{array}$ & $\begin{array}{cc}
 0 & 0 \\ a^{-1}bc
 &  0  \end{array}$ & $\begin{array}{cc} a &
b \\ 0 & c \end{array}$
\end{tabular} \right). \]
The group action on the torsion is then given by
\[ \begin{array}{lll}
A & = & A_{0}c^{-1}-C_{0}a^{-1}bc^{-1}-2ac^{-2}e \\
B & = & B_{0}ac^{-1} \\
C & = & C_{0}a^{-1} \\
D & = & D_{0}a^{3}c^{-3}+(1/2)B_{0}ab^{2}c^{-3}+(H_{0}-F_{0})a^{2}bc^{-3} \\
E & = & E_{0}a^{2}c^{-3} + F_{0}a^{3}c^{-4}e+G_{0}abc^{-3} \\
F & = & F_{0}a^{2}c^{-2} \\
G & = & G_{0}ac^{-2} \\
H & = & H_{0}a^{2}c^{-2}+B_{0}abc^{-2} \\
I & = & I_{0}c^{-1}+C_{0}a^{-1}bc^{-1} \\
J & = & J_{0}ac^{-2}-B_{0}a^{2}c^{-3}e-(1/2)C_{0}a^{-1}b^{2}c^{-2}-I_{0}bc^{-2}. 
\end{array} \]
We may then normalize
\[ A=0, \hspace{.1in} B=1, \hspace{.1in} C=1, \hspace{.1in} I=0. \]
It follows that
\[ v^{1}_{1} \equiv v^{2}_{2} \equiv v^{1}_{2} \equiv \beta \equiv 0 
\hspace{.22in} mod \hspace{.1in} base.  \]
We write 
\[ \begin{array}{lll}
v^{1}_{1} & = & A_{i}\eta^{i} \\
v^{2}_{2} & = & B_{i}\eta^{i} \\
v^{1}_{2} & = & C_{i}\eta^{i} \\
\beta & = & D_{i}\eta^{i}. \end{array} \]
Substituting these values back into the  equations (8) we obtain
\[ A_{5} = 2D_{4}+C_{4} \hspace{.5in} and \hspace{.5in} B_{4}=A_{1}+A_{4}.\]
The group has been reduced to the identity and so we have an e-structure.
The structure equations are \\
\\
$d\eta^{1} = (D_{1}+A_{2}-B_{2})\eta^{1}\eta^{2}+(A_{3}-B_{3})\eta^{1}\eta^{3}
-A_{1}\eta^{1}\eta^{4}+(A_{5}-B_{5})\eta^{1}\eta^{5}-D_{3}\eta^{2}\eta^{3} \\
\hspace*{.4in} +(1/2)(C_{4}-A_{5})\eta^{2}\eta^{4}-D_{5}\eta^{2}\eta^{5} $\\ \\
$d\eta^{2} = (1+2B_{1}-A_{1})\eta^{1}\eta^{2}+\eta^{1}\eta^{5}+(A_{3}-2B_{3})
\eta^{2}\eta^{3}-(A_{1}+B_{4})\eta^{2}\eta^{4}+(A_{5}-2B_{5})\eta^{2}\eta^{5}$ \\ \\
$d\eta^{3} = C_{1}\eta^{1}\eta^{2}+(1+B_{1})\eta^{1}\eta^{3} + \eta^{1}\eta^{4} +
(B_{2}-C_{3})\eta^{2}\eta^{3} -C_{4}\eta^{2}\eta^{4}-C_{5}\eta^{2}\eta^{5} \\
\hspace*{.4in} -(1+B_{4})\eta^{3}\eta^{4} -B_{5}\eta^{3}\eta^{5}$ \\ \\
$d\eta^{4} = I_{1}\eta^{1}\eta^{2}+I_{3}\eta^{1}\eta^{3}+A_{1}\eta^{1}\eta^{4}
+C_{1}\eta^{1}\eta^{5}-I_{2}\eta^{2}\eta^{3}+A_{2}\eta^{2}\eta^{4}
+C_{2}\eta^{2}\eta^{5} +A_{3}\eta^{3}\eta^{4}
\\ \hspace*{.4in} +(C_{3}-I_{4})\eta^{3}\eta^{5} 
+(C_{4}-A_{5})\eta^{4}\eta^{5}$ \\ \\
$d\eta^{5} = (C_{1}+I_{5})\eta^{1}\eta^{2}+\eta^{1}\eta^{4} +B_{1}\eta^{1}\eta^{5}
-C_{4}\eta^{2}\eta^{4} +(B_{2}-I_{6}-C_{5})\eta^{2}\eta^{5} +B_{3}\eta^{3}\eta^{5}
\\ \hspace*{.4in} +B_{4}\eta^{4}\eta^{5}-C_{3}\eta^{2}\eta^{3}$.
\\ \\ 
We obtain 24 local invariants.  

\end{document}